\documentclass[12pt, reqno]{amsart}
\usepackage[dvipsnames]{xcolor}
\usepackage{mathtools}
\mathtoolsset{showonlyrefs}
\usepackage{etoolbox}
\usepackage{amsthm, graphicx, microtype}
\usepackage[export]{adjustbox}

\usepackage{overpic}

\usepackage{blindtext}
\usepackage{latexsym}
\usepackage{txfonts}
\usepackage{lmodern}
\usepackage[pagewise]{lineno}
\usepackage[utf8]{inputenc}
\usepackage{exscale}
\usepackage{amsmath}
\usepackage{amssymb}
\usepackage{amsfonts}
\usepackage{mathrsfs}
\usepackage{amsbsy}
\usepackage{relsize}
\usepackage{comment}
\PassOptionsToPackage{reqno}{amsmath}
\usepackage{esint} 
\usepackage{amssymb} 
\usepackage{stmaryrd}
\usepackage{cite} 

\usepackage{tikz}

\usepackage{amsthm, graphicx, microtype}
\usepackage[export]{adjustbox}

\usepackage{overpic}
\usepackage[T2A,T1]{fontenc}

\definecolor{mypink1}{rgb}{0.858, 0.188, 0.478}
\definecolor{mypink2}{RGB}{219, 48, 122}
\definecolor{mypink3}{cmyk}{0, 0.7808, 0.4429, 0.1412}
\definecolor{mygray}{gray}{0.6}

\definecolor{gris75}{gray}{0.25}
\definecolor{violet}{rgb}{0.5,0,0.5}

\definecolor{BrickRed}{rgb}{0.58, 0.0, 0.83}

\definecolor{armygreen}{rgb}{0.29, 0.33, 0.13}
\definecolor{brass}{rgb}{0.71, 0.65, 0.26}
\definecolor{antiquefuchsia}{rgb}{0.57, 0.36, 0.51}
\definecolor{amethyst}{rgb}{0.6, 0.4, 0.8}
\definecolor{mauvetaupe}{rgb}{0.57, 0.37, 0.43}
\usepackage[colorlinks, linkcolor=purple, citecolor=blue, urlcolor=blue, hypertexnames=false]{hyperref}

\newcommand{\bS}{\mathbf{S}}
\newcommand{\bL}{\mathcal{A}}
\newcommand{\bk}{\mathbf{k}}
\newcommand{\bM}{\mathbf{M}}
\newcommand{\R}{\mathbb{R}}

\newcommand{\D}{\mathbb{D}}
\newcommand{\ts}{\mathsf{s}}

\newtheorem{prop}{Proposition}[section]
\newtheorem{thm}{Theorem}[section]
\newtheorem{rem}{Remark}[section]
\newtheorem{lem}{Lemma}[section]
\newtheorem{cor}{Corollary}[section]

\title[Subordinators and generalized heat kernels]{Subordinators and generalized heat kernels: Random time change and long time dynamics}

\author[N. Ajaber, A. Alshehri, H. Altamimi, M. Majdoub, E. Mliki]{N. Ajaber, A. Alshehri, H. Altamimi, M. Majdoub, E. Mliki}
\address{Department of Mathematics, College of Science, Imam Abdulrahman Bin Faisal University, P. O. Box 1982, Dammam, Saudi Arabia}
\address{Basic and Applied Scientific Research Center, Imam Abdulrahman Bin Faisal University, P.O. Box 1982, 31441, Dammam, Saudi Arabia}
\email{\sl \color{blue}{naljaber@iau.edu.sa}}
\email{\sl \color{blue}{asalshehri@iau.edu.sa}}
\email{\sl \color{blue}{haltamimi@iau.edu.sa}}
\email{\sl \color{blue}{mmajdoub@iau.edu.sa}}
\email{\sl \color{blue}{ermliki@iau.edu.sa}}

\begin{document}
\begin{abstract} 
This paper focuses on studying the long-time dynamics of the subordination process for a range of linear evolution equations, with a special emphasis on the fractional heat equation. By treating inverse subordinators as random time variables and employing the subordination principle to solve forward Kolmogorov equations, we explore the behavior of the solutions over extended periods. We provide a detailed description of the specific classes of subordinators suitable for conducting asymptotic analysis. Our findings not only extend existing research, but also enhance the results previously presented in \cite{KKS20, KKS21}.
\end{abstract}




\date{\today}

\maketitle

\section{Introduction}
This paper aims to investigate the long-term behavior of solutions for a specific class of linear evolution equations that emerge from random time changes utilizing various subordinators, such as {\tt gamma} and {\tt alpha} stable subordinators. We explore how these subordinators, denoted as $\bS = \left\{\bS_t\,; t \geq 0\right\}$, influence the dynamics of the solutions over extended time periods.

More specifically, our study focuses on examining the behavior of the subordinated solution, denoted as $v(x, t)$, for some evolution partial differential equations (PDEs) as time progresses. We analyze the large-time dynamics of the subordination process by convolving the solution with the density function $G_t$ associated with the subordinator $\bS$. This is accomplished through the analysis of the Cesàro mean of the subordinated solution. By examining the Cesàro mean, we gain insights into the long-term behavior and dynamics of the subordination process.

Understanding the behavior of the subordination process is of considerable significance in various physical models. Specifically, it has the potential to contribute to the development of valuable models that capture the concept of biological time in the evolution of species and ecological systems. By gaining insights into this process, we can enhance our understanding of complex biological phenomena and improve our ability to analyze ecological dynamics. This understanding can pave the way for advances in fields such as population dynamics, epidemiology, and ecological modeling, ultimately leading to more accurate predictions and informed decision making in real-world scenarios. For further explanations and detailed discussions on the topics covered, see \cite{MS2012, MK2000}.

Let us explain now the mathematical framework of the problem. Consider the linear evolution equation 
\begin{equation}\label{LEE}
\left \{ \begin{array}{rl}
\partial_t v(x, t)&=\bL\,v(x,t)\quad  \mbox{ in }\quad  \R^N\times (0,\infty),\\
v(x,0)&=\varphi(x) \quad\mbox{ in }\quad \mathbb{R}^N,
\end{array}\right.
\end{equation}
where $\bL$ is a Markov generator that acts on the functions $v(x,t)$. Assuming that $\eqref{LEE}$ has a solution $v$ such that $v(x,\cdot)\in L^1(\R_+)$, we introduce the subordination of $v(x,t)$ by a density function $G_t(\tau)$ as follows:
\begin{equation}
\label{v-E}
v^{E}(x,t):=\int_0^\infty\,v(x,\tau) G_t(\tau)\,d\tau, \quad x\in\R^N,\; t\geq 0.
\end{equation}
According to \cite{Baz}, we know that $v^{E}(x,t)$ is the solution of the general fractional differential equation
\begin{equation}
\label{FDE}
\left\{
\begin{matrix}
\left(\D^{(k)}_t v^{E}\right)(x, t)=\bL\,v^{E}(x,t)\\
v^{E}(x,0)= \varphi(x).
\end{matrix}
\right.
\end{equation}
Here $\D^{(k)}_t$ is the differential-convolution operator given by
\begin{equation}
\label{D-t-k}
\left(\D^{(k)}_t g\right)(t)=\frac{d}{dt}\int_0^t\, \mathbf{k}(t-s)\left(g(s)-g(0)\right)\,ds,
\end{equation}
where $\bk\in L^1_{loc}(\R_+)$ is a non-negative kernel. 

The main focus of this paper is to characterize the specific classes of subordinators that facilitate the derivation of time asymptotics for generalized fractional dynamics. In particular, our objective is to investigate the long-term behavior of $v^{E}(x,t)$, the solution of interest. To achieve this, we introduce the concept of Cesàro mean for $v^{E}(x,t)$:
\begin{equation}
\label{Cesàro}
\bM_t(v^{E}(x,t)):=\frac{1}{t}\int_0^t\, v^{E}(x,s)\,ds.
\end{equation}
Clearly, from \eqref{v-E} one can see that
\begin{equation}
\label{Cesàro-1}
\bM_t(v^{E}(x,t))=\int_0^\infty\, v(x,\tau) \bM_t\left(G_t(\tau)\right)\, d\tau.
\end{equation}
Our approach in this study is grounded in Laplace transform techniques, the Feller-Karamata-Tauberian theorem, and the integrability of the solution over time $v(x,t)$. This method can be applied to various models, including the fractional heat equation. In particular, our approach yields satisfactory results for dimensions $N\geq 2$. However, for the one-dimensional case, the obtained result is somewhat weaker. For further investigations and related studies, we recommend referring to \cite{KKS20, KKS21}.

We consider a broad class of admissible kernels, denoted $0 \leq \bk \in L^1_{\text{loc}}(\mathbb{R}_+)$, which are characterized in terms of their Laplace transforms $\mathcal{K}(\lambda)$ as follows.
\begin{equation}
    \label{Adm-1}
    \lambda \mathcal{K}(\lambda)\underset{\lambda\to  0^+}{\longrightarrow}\, 0,
\end{equation}
and, for some $\varrho\geq 0$, 
\begin{equation}
    \label{Adm-2}
    L(x):= x^{-\varrho}\,\mathcal{K}\left(x^{-1}\right)\quad\mbox{is a slowly varying function}\footnote{A function $L$ is said to be slowly varying function (SVF) if $\underset{x \to  \infty}{\lim} \frac{L(\lambda x)}{L(x)}=1, \;\; \mbox{for every}\; \lambda >0.$}.
\end{equation}
An example that is well known and commonly used is given by
\begin{equation}
    \label{k-mu}
    \bk(s)=\int_0^1\frac{s^{-\alpha} }{\Gamma(1-\alpha)}\mu(\alpha)d\alpha,
\end{equation}
where $\mu :[0,1]\to (0,\infty)$ is a continuous function. For further properties and asymptotic behavior of kernels defined by \eqref{k-mu}, we refer the reader to \cite{Koch}. In particular, \cite{Koch} provides detailed information on the Laplace transform of $\bk$ as follows:
\begin{equation}
\label{Lap-K}
\mathcal{K}(\lambda)=\int_0^1 \lambda^{\alpha-1} \mu(\alpha)\,d\alpha.
\end{equation}
In this case, it is worth noting that the slowly varying function $L$ is defined as follows:
\begin{equation}
    \label{L-SVF}
    L(x)=\int_0^1 x^{-\alpha} \mu(\alpha)d\alpha.
\end{equation}
Below, we provide a list of fundamental examples of kernel classes that will be utilized in our analysis.
 \begin{enumerate}
\item[{\bf 1})] The  $\theta$-stable subordinator class (${\bf C}_1$):
\begin{equation*}
\mathcal{K}(\lambda) =\lambda^{\theta -1}, \quad 0<\theta<1.\end{equation*}
\item[{\bf 2})]  The distributed order derivatives class (${\bf C}_2$):
\begin{equation*}
\mathcal{K}(\lambda)\; \underset{\lambda\to  0^+}{\sim}\; C\lambda^{ -1}\left(-\ln\lambda\right)^{-\kappa}, \quad C,\kappa >0.\end{equation*}
\item[{\bf 3})]  The  inverse Gamma subordinator class (${\bf C}_3$):
\begin{equation*}
\mathcal{K}(\lambda) =\frac{\sqrt{b}}{\lambda} \left(2\sqrt{2\lambda+a}-\sqrt{a}\right), \quad a\geq 0,\;\;b>0.\end{equation*}
\item[{\bf 4})] The Gamma subordinator class (${\bf C}_4$):
\begin{equation*} \mathcal{K}(\lambda) = \frac{a}{\lambda}\ln\left(1+\frac{\lambda}{b}\right), \quad a, b >0.\end{equation*}
\item[{\bf 5})] The tempered stable subordinator class (${\bf C}_5$):
\begin{equation*} \mathcal{K}(\lambda) = \frac{(\lambda+\beta)^\theta-\beta^\theta}{\lambda}, \quad \beta >0.\end{equation*}
\end{enumerate}
\begin{rem}\label{Classes}
    ~{\rm 
    \begin{itemize}
        \item[(i)] All classes (${\bf C}_i$), $1\leq i\leq 5$, satisfy assumption \eqref{Adm-1} except for (${\bf C}_3$) with $a>0$. It is worth noting that when $a=0$ in (${\bf C}_3$), we obtain the class (${\bf C}_1$) as a special case with $\theta=1/2$.
        \item[(ii)] The hypothesis \eqref{Adm-2} holds for all classes (${\bf C}_i$), $1\leq i\leq 5$, with specific values of $\varrho$. For (${\bf C}_1$), we have $\varrho=1-\theta$ and $L(x)=1$, for (${\bf C}_2$), we have $\varrho=1$ and $L(x)=C\left(\ln x\right)^{-\kappa}$, for (${\bf C}_3$) we have $\varrho=1$ and $L(x)=\sqrt{b}\left(2\sqrt{2/x +a}-\sqrt{a}\right)$, for (${\bf C}_4$), we have $\varrho=0$ and $L(x)=ax\ln\left(1+\frac{1}{bx}\right)$ and for (${\bf C}_5$), we have $\varrho=0$ and $L(x)=x\left(\left(1/x+\beta\right)^\theta-\beta^\theta\right)$.
    \end{itemize} 
    }
\end{rem}

Henceforth, we will consistently represent the Lebesgue norm in $L^p$, applicable to both space and time variables, as $\|\cdot\|_p$ for $1\leq p\leq \infty$.
The following statement presents our first main result:
\begin{thm}
\label{Main}
Suppose that the Markov generator $\mathcal{A}$ satisfies the condition that for any $0 \leq \varphi \in L^1(\mathbb{R}^N)$, the solution $v$ of \eqref{LEE} is non-negative, and $v(x, \cdot) \in L^1(0, \infty)$.
Moreover, if we assume that the conditions \eqref{Adm-1} and \eqref{Adm-2} are fulfilled, then:
\begin{equation}
    \label{main-res}
    \bM_t(v^{E}(x,t))\, \underset{t\to \infty }{\backsim}\, \frac{\|v(x)\|_{1}}{\Gamma(\varrho+1)} \,\bigg(t^{-1} \mathcal{K}\left(t^{-1}\right)\bigg).
\end{equation}
Here, $v^{E}(x,t)$ represents the subordination of the solution $v(x,t)$ as defined in \eqref{v-E}.
\end{thm}
\begin{rem}
~{\rm 
\begin{itemize}
    \item[(i)] The asymptotic behavior result, stated as \eqref{main-res}, can be alternatively expressed as
    \begin{equation}
    \label{main-res-equivalent}
    \bM_t(v^{E}(x,t))\, \underset{t\to \infty }{\backsim}\, \frac{\|v(x)\|_{1}}{\Gamma(\varrho+1)} \,\left(\int_0^1 t^{-\alpha} \mu(\alpha)d\alpha\right),
\end{equation}
where the function $\mu$ is given by \eqref{k-mu}.
\item[(ii)] If we assume $ \lambda \mathcal{K}(\lambda)\underset{\lambda\to 0^+}{\longrightarrow}\,\ell>0$ instead of \eqref{Adm-1}, then the requirement $v(x,\cdot)\in L^1(0,\infty)$ can be omitted. In that case, the conclusion \eqref{main-res} remains valid, provided that $\|v(x)\|_{1}$ is replaced with
$$
\|v(x)\|_{1,\ell}:=\int_0^\infty\,{\rm e}^{-\ell\tau}\,v(x,\tau)\,d\tau.
$$
    \item[(iii)] Our primary contribution is encapsulated in Theorem \ref{Main}, where we employ a unified approach to investigate a broad class of operators and kernels.
\end{itemize}
}  
\end{rem}
For specific choices of kernel classes, Theorem \ref{Main} provides the following asymptotic estimates.
\begin{cor}
    \label{Spec-Kern}
    Under the assumptions of Theorem \ref{Main} {\rm(}except for class {\rm(}$\mathbf{C}_3${\rm)} where the assumption $v(x,\cdot)\in L^1$ is omitted{\rm)}, the long-time behavior of the Cesàro mean of $v^{E}(x,t)$ as $t\to\infty$ is expressed as:
    \begin{equation}
    \label{Calss-Examp}
    \begin{split}
\bM_t(v^{E}(x,t))\, \underset{t\to \infty }{\backsim}\, \left\{\begin{array}{llll}   \frac{\|v(x)\|_1}{\Gamma(2-\theta)}\,t^{-\theta} \quad &\text{if}\quad \mathbf{k}\in\; (\mathbf{C}_1)\\
\mathtt{C}\,\|v(x)\|_1\,\left(\ln t\right)^{-\kappa} \quad &\text{if}\quad \mathbf{k}\in\; (\mathbf{C}_2)\\
 2\sqrt{ab}\,\|v(x)\|_{1,\sqrt{ab}}  \quad &\text{if}\quad \mathbf{k}\in\; (\mathbf{C}_3)\\
 \frac{a\|v(x)\|_1}{b}\, t^{-1} \quad &\text{if}\quad \mathbf{k}\in\; (\mathbf{C}_4)\\
 \theta\beta^{\theta-1}\,\|v(x)\|_1\, t^{-1} \quad &\text{if}\quad \mathbf{k}\in\; (\mathbf{C}_5).
\end{array}\right.
\end{split}
\end{equation}
\end{cor}
Theorem \ref{Main} leads to the following consequential result.  
\begin{thm}
\label{frac-lap}
Consider the operator $\mathcal{A}=-(-\Delta)^{\ts}$ where  either {\rm(}$N\geq 2$ and $0<\mathtt{s}<1${\rm)} or {\rm(}$N=1$ and {$0< \ts<1/2$}{\rm)}. If $0 \leq \varphi \in L^1\cap L^\infty(\mathbb{R}^N)$, the solution $v$ of \eqref{LEE} is non-negative, and $v(x, \cdot) \in L^1(0, \infty)$. Consequently, the asymptotic behavior described in \eqref{main-res} holds true, given the conditions \eqref{Adm-1} and \eqref{Adm-2} imposed on the kernel.
\end{thm}
\begin{rem}
~{\rm 
\begin{itemize}
    \item[(i)] The case $\ts=1$ was explored in \cite{KKS20, KKS21} across all dimensions.
    \item[(ii)] The non-negativity of the solution $v$ only necessitates $0<\ts\leq 1$ (see Lemma \ref{Positive} below). 
    \item[(iii)] In order to guarantee that $v(x,\cdot)\in L^1(1,\infty)$, we need to impose an additional condition: either $N\geq 2$ or, in the case of $N=1$, we require $0<\ts<1/2$ along with $\varphi\in L^1$. See Proposition \ref{Smooth-eff} below.
    \item[(iv)] The requirement of $\varphi\in L^\infty$ is used to guarantee the integrability of $v(x,\cdot)$ near $t=0$.
    \item[(v)] In addition to the nonlocal operator discussed in Theorem \ref{frac-lap}, another class of nonlocal operators that can be addressed using the same framework is represented by $\mathcal{A}v= \mathbf{J}\ast v-v$. Here, $\mathbf{J}$ denotes a radial probability density function that satisfies certain additional assumptions. For further details and insights into this class of nonlocal diffusion problems, refer to \cite{NLDP}.
\end{itemize}
}  
\end{rem}
When considering the scenario where $N=1$ and $1/2\leq \ts<1$, we obtain a result that is slightly weaker, as presented below.
\begin{thm}
    \label{N=1}
Assuming $N=1$ and $1/2\leq \ts<1$, and considering a kernel $\mathbf{k}$ that satisfies assumptions \eqref{Adm-1} and \eqref{Adm-2}, we have the following  long-term behavior of the Cesàro mean of $v^{E}(x,t)$:
\begin{equation}
    \label{LTB-N=1}
   \bM_t(v^{E}(x,t))\, \asymp \, \frac{t^{\frac{\varrho-1}{2\ts}}}{\Gamma(1+\varrho)}\, \left(L(t)\right)^{\frac{1}{2\ts}}.
\end{equation}
Here, the notation $\asymp$ signifies that the ratio under consideration is bounded both from above and below. 
\end{thm}

The structure of the paper is as follows. In Section \ref{S2}, we provide a recap of the essential background information required for the subsequent proofs. Section \ref{S3} is dedicated to presenting the proofs of our main results, namely Theorem \ref{Main}, Corollary \ref{Spec-Kern}, Theorem \ref{frac-lap}, and Theorem \ref{N=1}.

\section{Preliminaries}
\label{S2}

\subsection{Notations}
\begin{itemize}
    \item[$\triangleright$] We write $X\lesssim Y$ or $Y\gtrsim X$ to denote the estimate $X\leq CY$ for some constant $C>0$. 
    \item[$\triangleright$] The notation $f(z)\underset{z\to z_0}{\backsim}g(z)$ means that $\displaystyle\lim_{z\to z_0}\frac{f(z)}{g(z)}=1$, where $-\infty\leq z_0\leq \infty$.
    \item[$\triangleright$] The notation $f(z)\;\asymp\;g(z)$ indicates that the ratio is bounded by a constant factor both from above and below.
\item[$\triangleright$] The Laplace transform of a function $f: [0,\infty)\to\R$ is defined as follows
   $$
    \mathcal{L}(f)(\lambda):=\int_0^\infty\, {\rm e}^{-\lambda t}\,f(t)\,dt,\;\; \lambda>0.
    $$
    \item[$\triangleright$]  We will use the notation $\R_+=(0,\infty)$ to represent the set of nonnegative real numbers.
     \item[$\triangleright$] To simplify notation, we abbreviate the Lebesgue norm on $\R_+$ as $\|\cdot\|_{1}:= \|\cdot\|_{L^1(0,\infty)}$

\item[$\triangleright$] Throughout the entire paper, we will use the letter $C$ to represent various positive constants that are not significant in our analysis and may vary from line to line.
\end{itemize}

\subsection{Subordinators}
{A subordinator is a process with stationary and independent non-negative increments starting at zero.
Subordinators are a special class of L\'{e}vy  processes taking values in $[0, \infty)$ and their sample
paths are non-decreasing; this is a type of stochastic process that is used to model random phenomena that have jumps (see  \cite{Bern,Ber} for more details).  Let $\bS = \left\{\bS_t,\; t \geq 0\right\}$ be a subordinator. } The infinite divisibility of the law of $\bS$ implies that
its Laplace transform can be expressed in the form
\begin{eqnarray*}
E(e^{-\lambda \bS_{t}})=e^{-t \Phi(\lambda)}=e^{-t \lambda  \mathcal{K}(\lambda)}, \quad \lambda\geq0,
\end{eqnarray*}
where $\Phi : [0, \infty) \to [0, \infty)$, called the Laplace exponent, is
a Bernstein function \cite{Bernstein}. 
Such a Laplace exponent $\Phi$ can be expressed as
\begin{eqnarray*}
\Phi(\lambda)=\int_{(0, \infty)} (1-e^{-\lambda \tau }) \, \sigma(d\tau) <\infty,
\end{eqnarray*}
which is known as the L\'{e}vy-Khintchine formula for the subordinator $\bS$. Where  $\sigma$ are a measure on the positive real half-line such that 
\begin{eqnarray*}
\int_{(0, \infty)} (1 \land \tau ) \, \sigma(d\tau) <\infty.
\end{eqnarray*}
The kernel $\bk$ is related to the subordinator $\bS$ via the L\'{e}vy measure $\sigma$, namely if we set 
$$\bk(t)=\sigma((t, \infty)), \quad t\geq  0.$$
For any $\lambda \geq 0$ we have
$$\int_0^\infty e^{-\lambda t}\int_0^t d\sigma(s)\,ds=\int_0^\infty \int_0^s e^{-\lambda t} dt\,d\sigma(s)=\frac{1}{\lambda}\Phi(\lambda)=\mathcal{K}(\lambda).$$
Denote by $E$ the inverse process of the subordinator $\bS$, 
$$
E_t:=\mbox{inf}\{ s\geq0;\, \bS_t\geq t\}.
$$
The marginal density of $E_t$ is the function $G_t(\tau)$, more precisely
$$G_t(\tau)\, d\tau=\partial_\tau P(E_t \leq \tau )=\partial P(\bS \geq t).$$
\subsubsection{ $\theta$-stable subordinator}
 Let $\bk$ be the kernal corresponding to the Caputo-Djrbashian fractional derivative $\D^{(\theta)}_t$ of order $\theta \in (0, 1).$ Then its Laplace transform is given by 
\begin{equation*} \mathcal{K}(\lambda)=\frac{1}{\Gamma(1-\theta)}\int_0^\infty\, {\rm e}^{-\lambda t}\,t^{-\theta}\,dt= \lambda^{\theta-1},
    \end{equation*}
where $\bk(t)=\frac{t^{-\theta}}{\Gamma(1-\theta)}.$
It is easy to verify that (\ref{Adm-1}) holds for $\mathcal{K}.$

\subsubsection{Gamma Subordinator}
 Let $\bk$ be the kernal given by 
\begin{equation*} \bk(t)=a\Gamma(0, bt), \quad a,b>0,\end{equation*}
    where $\Gamma(\nu, x):\displaystyle{\int_x^\infty t^{\nu-1}e^{-t}\,dt}$ is the upper incomplete Gamma function. 
The Laplace transform of $\bk$ is given by 
\begin{equation*} \mathcal{K}(\lambda) = \frac{a}{\lambda}\ln\left(1+\frac{\lambda}{b}\right), \;\; a, b,\lambda >0.\end{equation*}

It is easy to see that (\ref{Adm-1}) holds for $\mathcal{K}.$

\subsubsection{Inverse Gamma Subordinator} Let $a\geq 0$ and $b>0$ be given and definite the kernel $\bk$ by 
\begin{equation*} \bk(t)=\sqrt{\frac{b}{2\pi}}\left(\frac{2}{\sqrt{t}}e^{-\frac{at}{2}}-\sqrt{2a\pi}(1-\mbox{erf}(z))\right),\quad z:=\sqrt{\frac{at}{2}},\end{equation*}
where $\mbox{erf}(z):=\frac{2}{\sqrt{\pi}}\displaystyle\int_0^z e^{-t^2}\,dt$ is the error function. 
The Laplace transform of $\bk$ is given by 
\begin{equation*}
\mathcal{K}(\lambda) =\frac{\sqrt{b}}{\lambda} \left(2\sqrt{2\lambda+a}-\sqrt{a}\right), \quad a\geq 0,\;\;b, \lambda>0.\end{equation*}

The assumption (\ref{Adm-1}) holds for $\mathcal{K}$ when the value $0$ in \eqref{Adm-1} is replaced with $\sqrt{ab}$.

\subsubsection{Tempered Stable Subordinator}
 Let us consider a tempered stable subordinator $D_{\beta}(t)$ with index $\theta\in (0,1)$. The density function of $D_{\beta}(t)$ can be expressed as:
 $$
 f_\beta(x,t)={\rm e}^{-\beta\,x+\beta^{\theta}t}\,f(x,t),\quad \beta>0,
 $$
 where $\beta > 0$ and $f(x,t)$ represents the density of the $\theta$-stable subordinator. It is worth noting that $f_\beta(x,t)$ possesses finite moments and exhibits infinite divisibility. However, it is not self-similar.

Moreover, we can obtain the Laplace transform of the density function $f_\beta(x,t)$, given by:
$$
\mathcal{L}_x\left(f_\beta(x,t)\right)={\rm e}^{-t\left((\lambda+\beta)^{\theta}-\beta^{\theta}\right)}={\rm e}^{-t\lambda\mathcal{K}(\lambda)}.
$$
For further details and explanations, refer to \cite{KGW, Mee2013}.
\subsection{Fractional Heat Semigroup}
Consider the linear homogeneous fractional heat equation 
\begin{equation}\label{HLE}
\partial_t v+(-\Delta)^{\ts}v=0,\quad v(0)=\varphi.
\end{equation}
It is well known that the operator $(-\Delta)^{\ts}$, $\ts>0$ serves as the generator of a semi-group $e^{-t(-\Delta)^{\ts}}$, whose kernel $\mathscr{E}_{\ts}$ is smooth, radial, and adheres to the scaling property:
\begin{equation}\label{kernel}
\mathscr{E}_{\ts}(x,t)=t^{-\frac{N}{2\ts}}\mathscr{E}_\ts(t^{-\frac{1}{2\ts}}x,1):=t^{-\frac{N}{2\ts}} \mathscr{K}_\ts\left(t^{-\frac{1}{2\ts}}x\right),
\end{equation}
where 
\begin{equation}
    \label{K-s}
    \mathscr{K}_\ts(x)=(2\pi)^{-N/2}\int_{\R^N}\,e^{i x\cdot\xi}\,e^{-|\xi|^{2\ts}}\,d\xi.
\end{equation}
Therefore, the solution to \eqref{HLE} can be formally obtained through convolution as
\begin{equation}
\label{convol}
  v(x,t)=\left(e^{-t(-\Delta)^{\ts}}\right)\varphi(x)=\bigg(t^{-\frac{N}{2\ts}} \mathscr{K}_\ts\left(t^{-\frac{1}{2\ts}}\cdot\right)\ast \varphi\bigg)(x), 
\end{equation}
 whenever this representation is meaningful. The explicit expression of the kernel $\mathscr{E}_{\ts}$ is known for particular cases $\ts=1$ and $\ts=1/2$. When $\ts=1$ we obtain the standard Gaussian kernel:
 \begin{equation}
     \label{Gauss}
  \mathscr{E}_{1}(x,t)=(4\pi t)^{-N/2}\,{\rm e}^{-\frac{|x|^2}{4t}},   
 \end{equation}
 with profile $\mathscr{K}_{1}(x)=(4\pi)^{-N/2}\,{\rm e}^{-\frac{|x|^2}{4}}.$ For $\ts=1/2$, we get the well known Poisson kernel
 \begin{equation}
     \label{Poisson}
   \mathscr{E}_{1/2}(x,t)=  \frac{\Gamma(N+1/{2})\,t}{\pi^{\frac{2N+1}{2}}\big(t^2+|x|^2\big)^{\frac{N+1}{2}}},
 \end{equation}
with profile $\mathscr{K}_{1/2}(x)=\frac{\Gamma(N+1/{2})}{\pi^{\frac{2N+1}{2}}\big(1+|x|^2\big)^{\frac{N+1}{2}}}$. Here, the symbol $\Gamma$ stands for the Gamma function. 

While the explicit form of $\mathscr{K}_{\ts}$ is not known for every $\ts \in (0,1)$, we do have a positivity property, as stated below.
\begin{lem}
    \label{Positive}
Let $N\geq 1$ and $\ts\in(0,1)$. Then for all $x\in\R^N$, we have 
\begin{equation}
         \label{kthetaest}
        {c_{\ts}}\,\left(1+|x|\right)^{-N-2\ts}\,\leq\,\mathscr{K}_{\ts}(x)\leq  C_{\ts}\,\left(1+|x|\right)^{-N-2\ts},
     \end{equation}
     where  $c_{\ts},\, C_{\ts}$ are positive constant depending only on $N$ and $\ts$.
As a result, we get $\mathscr{K}_{\ts}\in L^p(\R^N)$ for any $1\leq p\leq \infty$.
\end{lem}
The proof of Lemma \ref{Positive} can be found, for instance, in \cite[p. 395]{Alonso2021}. See also \cite[Theorem 2.1, p. 263]{BG1960}.  It's worth noting that the crux of the proof relies on the formula
\begin{equation}
    \label{Posit-form}
    \mathscr{K}_{\ts}(x)=\displaystyle\int_0^\infty
    \,(4\pi t)^{-N/2}\,{\rm e}^{-\frac{|x|^2}{4t}} \chi(t)\,dt,
\end{equation}
where $\chi(t)$ represents the inverse Laplace transform of  $0<\lambda\mapsto{\rm e}^{-\lambda^{\ts}}$. It is noteworthy that the positivity of the kernel $\mathscr{K}_{\ts}$ was asserted in \cite[p. 263]{BG1960} without providing a proof.
\begin{rem}
    {\rm  As indicated by \cite[Theorem 2.1]{BG1960}, the constant $c_{\ts}$ in the lower bound of \eqref{kthetaest} approaches $0$ as $\ts$ approaches $1$. Specifically, the lower bound estimate in \eqref{kthetaest} becomes invalid when $\ts=1$ due to the exponential decay of $\mathscr{E}_{\ts}(\cdot, t)$.}
\end{rem}

By using \eqref{kernel}, \eqref{convol}, and \eqref{kthetaest}, we can easily establish both a lower and an upper bound for $v(x,t)$ as follows.
\begin{prop}
    \label{bound-v}
    Let $N\geq 1$, $\ts\in (0,1)$ and $0\leq \varphi\in L^1(\R^N)$. For any compact set $K\subset \R^N$, there exists a constant $C=C(K,\varphi,N)>1$, such that
    \begin{equation}
        \label{v-bound}
        \frac{1}{C} t^{-\frac{N}{2\ts}}\,\leq\, v(x,t)\,\leq\,{C} t^{-\frac{N}{2\ts}},\quad t\geq 1,\,\, x\in K.
    \end{equation}
\end{prop}

By utilizing \eqref{convol}, Lemma \ref{Positive}, and Young's inequality, we readily derive the subsequent $L^p-L^q$ estimate.
\begin{prop}\label{Smooth-eff}
For $N\geq 1$ and $0<\ts<1$, there exists a positive constant $C=C(N,\ts)$ such that for all $1\leq p\leq q\leq \infty$, the following inequality holds:
\begin{equation}\label{Smooth-est}
\|e^{-t(-\Delta)^{\ts}}\varphi\|_{q}\leq C\,t^{-\frac{N}{2\ts}\left(\frac{1}{p}-\frac{1}{q}\right)}\,\|\varphi\|_{p},\quad t>0,\;\; \varphi\in L^p.
\end{equation}
\end{prop}
 By exploiting the property that $e^{-t(-\Delta)^{\ts}}$ forms a $C_0$-semigroup, we derive a result analogous to the one established in \cite[Lemma 2.1]{Laister} for $\ts=1$.
 \begin{prop}
    \label{C-0}
    Suppose $0<\ts<1$, $1\leq q<r\leq\infty$, and let $\delta=\frac{N}{2\ts}\left(\frac{1}{q}-\frac{1}{r}\right)$. If $\varphi\in L^q(\R^N)$, then the following holds:
    \begin{equation}
        \label{C-0-0}
        \lim_{t\to 0}\, t^{\delta}\,\left\|e^{-t(-\Delta)^{\ts}}\varphi\right\|_{L^r}=0.
    \end{equation}
 \end{prop}

The asymptotic behavior of the solution $v$ to \eqref{HLE} was explored in \cite{Vaz}. The principal result in \cite{Vaz} demonstrates that any integrable solution of \eqref{HLE} exhibits behavior similar to the fundamental solution $\mathscr{E}_{\ts}(x,t)$ provided in \eqref{kernel}.
\begin{thm}{\cite[Theorem 3.1, p. 1220]{Vaz}}
    \label{thm-Vaz}
Consider $\ts \in (0,1)$ and $\varphi \in L^1(\mathbb{R}^N)$ with $\mathbf{m}:=\int_{\mathbb{R}^N}\varphi(x)dx\neq 0$. Then, the solution $v(x,t)$ to \eqref{HLE} satisfies the following asymptotic estimates:
\begin{equation}
    \label{Asymp-1}
    \Big\|v(t)-\mathbf{m}\mathscr{E}_{\ts}(t)\Big\|_{1}\, \underset{t \to  \infty}{\longrightarrow}\, 0
\end{equation}
and
\begin{equation}
    \label{Asymp-2}
    t^{N/2\ts}\,\Big\|v(t)-\mathbf{m}\mathscr{E}_{\ts}(t)\Big\|_{\infty}\, \underset{t \to  \infty}{\longrightarrow}\, 0.
\end{equation}
\end{thm}
As emphasized in \cite{Vaz}, it is not necessary for $v$ to be non-negative, and consequently, the mass $\mathbf{m}$ may have any sign.
\subsection{Karamata's Tauberian theorem}
The Feller-Karamata Tauberian theorem is a powerful result in mathematical analysis that establishes a connection between the behavior of a function and its Laplace transform. See \cite[Sec. 1.7]{Bingham} and \cite[XIII, Sec. 1.5]{Feller} for a general statement. Below, we present a simplified version of the Feller-Karamata Tauberian theorem that adequately serves our purpose.
\begin{thm}
\label{FKT}
Let $U:[0,\infty)\to\R$ be a monotone non-decreasing right-continuous function satisfying $\mathcal{L}(U)(\lambda)<\infty$ for all $\lambda>0$. Then, the following assertions are equivalent:
\begin{equation}
\label{U-t}
U(t)\sim \frac{C t^\rho}{\Gamma(\rho+1)}L(t)\quad\text{as}\quad t\to\infty,
\end{equation}
\begin{equation}
\label{w-lambda}
\mathcal{L}(U)(\lambda)\sim C \lambda^{-1-\rho}L\left(\frac{1}{\lambda}\right)\quad\text{as}\quad \lambda\to 0^+.
\end{equation}
Here, $C>0$, $\rho\geq 0$, and $L$ represents a slowly varying function.
\end{thm}
\begin{rem}
~
\begin{itemize}
\item[(i)] {\rm { In essence, Theorem \ref{FKT} establishes a direct connection between the asymptotic tendencies of $U(t)$ as $t$ approaches infinity and the corresponding behavior of its Laplace transform as $\lambda$ tends to $0$. The constants $C$, $\rho$, and the slowly varying function $L$ capture the precise quantitative relationship between the two asymptotic behaviors.}}
    \item[(ii)] {\rm {The theorem's validity can be confirmed for the function $U(t)=\frac{t^\rho}{\Gamma(\rho+1)}$, where $\rho\geq 0$, resulting in  $L(t)=1$ a clearly slowly varying function. By choosing $U(t)=\frac{t^\rho}{\Gamma(\rho+1)}$ it becomes evident that $\mathcal{L}(U)(\lambda)=\lambda^{-1-\rho}$. }}
    \end{itemize}
\end{rem}

The following asymptotic behavior of integral uses the incomplete Gamma function and it is needed in the proof of Theorem \ref{N=1}.
\begin{lem}
    \label{IGF}
    Consider $0<\theta\leq 1$, and define
    $$
    \mathbf{F}(\varepsilon)=\int_1^\infty\,\tau^{-\theta}{\rm e}^{-\varepsilon\tau}\, d\tau,\quad \varepsilon>0.
    $$
    Then, the following asymptotic holds 
    \begin{equation}
        \label{IGF-asymp}
         \mathbf{F}(\varepsilon)\,\sim \,\varepsilon^{\theta-1}\, \Gamma(1-\theta)\quad \text{as}\quad \varepsilon\to 0.
    \end{equation}
\end{lem}

\section{Proof of the main result}
\label{S3}
This section is devoted to the proof of Theorem \ref{Main}. First, recall a crucial relation between the Laplace transform of the kernel $k$
 and the density $G_t(\tau)$.
 \begin{lem}
    \label{G-K}
    The $t-$Laplace transform of $G_t(\tau)$ is given by
    \begin{equation}
\label{G-K-lambda}
        \int_0^\infty {\rm e}^{-\lambda t}G_t(\tau)\,dt =\mathcal{K}(\lambda){\rm e}^{-\tau \lambda\mathcal{K}(\lambda)},
    \end{equation}
   where $\mathcal{K}(\lambda)$ is the Laplace transform of the kernel $k$.
 \end{lem}
 The proof of the above Lemma can be found in \cite[Lemma 2.2]{KKS20} for instance.
 
Let $0 \leq \varphi \in L^1(\mathbb{R}^N)$, and  $v$ the non-negative solution of \eqref{LEE}. It follows that the subordination $v^{E}$ of $v$ by the density  $G_t(\tau)$ given by \eqref{v-E} is also non-negative. Let us define 
\begin{equation}
\label{U}
\mathbf{U}(t)=\displaystyle\int_0^t\,v^{E}(x,s)\,ds=t \bM_t(v^{E}(x,s)).\end{equation} 
Here, $\bM_t$ is the Cesàro mean given by \eqref{Cesàro-1}.  
By using the formula \eqref{G-K-lambda}, we obtain that 
\begin{equation}
\label{w-K}
\frac{\lambda\,w(\lambda)}{\mathcal{K}(\lambda)}= \int_0^\infty {\rm e}^{-\tau \lambda \mathcal{K}(\lambda) } v(x, \tau)\,d\tau,
\end{equation}
where $w(\lambda):=\mathcal{L}(\mathbf{U})(\lambda)$ denotes the Laplace transform of $\mathbf{U}$.

Thanks to \eqref{Adm-1} and the fact that $v(x,\cdot)\in L^1(0,\infty)$, the Dominate Convergence Theorem yields 
\begin{equation}
\label{lambda -0}
\underset{\lambda \to  0^+}{\lim} \left(\frac{\lambda w(\lambda)}{\mathcal{K}(\lambda)}\right)=\int_0^\infty v(x, \tau) \, d\tau=\|v(x)\|_{1}.
\end{equation}
Hence,
\begin{equation}
    \label{w-lam-0}
    w(\lambda)\, \underset{\lambda \to  0^+}{\sim}\,\|v(x)\|_{1}\,\frac{\mathcal{K}(\lambda)}{\lambda}.
\end{equation}
Referring to \eqref{Adm-2}, the above asymptotic formula \eqref{w-lam-0} can be expressed alternatively as
\begin{equation}
    \label{w-lam-L}
w(\lambda)\, \underset{\lambda \to  0^+}{\sim}\,\|v(x)\|_{1}\,\lambda^{-1-\varrho}\,L\left(\frac{1}{\lambda}\right).
\end{equation}
Here, as stated in \eqref{Adm-2}, $\varrho\geq 0$ and $L$ is a slowly varying function (SVF) at infinity. 
Due to the non-negativity of $v^{E}$, we see that $\mathbf{U}$ is a continuous non-decreasing function. Hence, we  can apply Feller-Karamata-Tauberian theorem and obtain that 
\begin{equation}
\label{FKT-App}
\bM_t(v^E(x, t))\, \underset{t \to \infty}{\sim} \,\frac{\|v(x)\|_{1}}{\Gamma(\varrho +1)} t^{\varrho-1}L(t). 
\end{equation}
From the hypotheses \eqref{Adm-2} we easily derive \eqref{main-res}. This finishes the proof of Theorem \ref{Main}.
\subsection{Proof of Corollary \ref{Spec-Kern}}
Remark \ref{Classes} highlights that all the classes under consideration satisfy assumptions \eqref{Adm-1} and \eqref{Adm-2}, except for (${\bf C}_3$) where the value $0$ in \eqref{Adm-1} is replaced with $\sqrt{ab}$. Consequently, the desired conclusion \eqref{Calss-Examp} directly follows from \eqref{main-res}. This concludes the proof of Corollary \ref{Spec-Kern}.
\subsection{Proof of Theorem \ref{frac-lap}}
Since $\varphi\geq 0$, Lemma \ref{Positive} implies that $v\geq 0$. Moreover, due to the fact that $\varphi\in L^1$, it can be readily inferred from \eqref{kernel} and \eqref{kthetaest} that
\begin{equation}
\label{v-L1}
0\leq v(x,t)\leq C_{\ts}\,\|\varphi\|_1\,t^{-\frac{N}{2\ts}},\quad t>0,\,\,\, x\in\R^N.
\end{equation}
Consequently, we can conclude that $v(x,\cdot)\in L^1(1,\infty)$ for $N>2\ts$. This condition is satisfied when $N\geq 2$ or $N=1$ and $\ts<1/2$.

Next, by employing the fact that $\varphi\in L^\infty$ in conjunction with \eqref{Smooth-eff}, we obtain
\begin{equation}
\label{v-Linfty}
0\leq v(x,t)\leq \|\varphi\|_\infty,\quad t>0,\,\,x\in\R^N.
\end{equation}
From \eqref{v-Linfty}, it is evident that $v(x,\cdot)\in L^1(0,1)$. Consequently, the solution $v$ is nonnegative and satisfies $v(x,\cdot)\in L^1(0,\infty)$ whenever $N\geq 2$ or $N=1$ and $\ts<1/2$. Thus, the proof of Theorem \ref{frac-lap} can be concluded by invoking Theorem \ref{Main}.

\subsection{Proof of Theorem \ref{N=1}}
The proof of Theorem \ref{N=1} relies on the utilization of the incomplete Gamma function (Lemma \ref{IGF}) and Proposition \ref{v-bound}. By employing Lemma \ref{IGF}, Proposition \ref{v-bound}, and decomposing the integral as follows:
$$
\int_0^\infty {\rm e}^{-\tau \lambda \mathcal{K}(\lambda) } v(x, \tau)\,d\tau=\int_0^1 {\rm e}^{-\tau \lambda \mathcal{K}(\lambda) } v(x, \tau)\,d\tau+\int_1^\infty {\rm e}^{-\tau \lambda \mathcal{K}(\lambda) } v(x, \tau)\,d\tau,
$$
we deduce that 
\begin{equation}
    \label{N-1}
    \int_0^\infty {\rm e}^{-\tau \lambda \mathcal{K}(\lambda) } v(x, \tau)\,d\tau\,\asymp\, \left( \lambda \mathcal{K}(\lambda)\right)^{\frac{1-2\ts}{2\ts}}\Gamma\left(\frac{2\ts-1}{2\ts}\right),
\end{equation}
uniformly on compact subsets of $\R^N$. By invoking \eqref{w-K}, we obtain
\begin{equation}
    \label{N-1-1}
    w(\lambda)\,\asymp\, \lambda^{\frac{1}{2\ts}-2}\, \left(\mathcal{K}(\lambda)\right)^{\frac{1}{2\ts}}.
\end{equation}
Now, using \eqref{Adm-2}, the above asymptotic can be expressed as
\begin{equation}
    \label{N-1-2}
    w(\lambda)\,\asymp\,\lambda^{-1-(1+\frac{\varrho-1}{2\ts})}\, \left(L({\frac{1}{\lambda})}\right)^{\frac{1}{2\ts}},
    \end{equation}
    where $L$ is given by \eqref{Adm-2}. Hence, by applying Karamata's Tauberian theorem we obtain
    \begin{equation}
    \label{N-1-3}
\mathbf{U}(t) \asymp \frac{t^{1+\frac{\varrho-1}{2s}}}{\Gamma(1+\varrho)}\, \left(L(t)\right)^{\frac{1}{2\ts}}.
    \end{equation}
The proof of Theorem \ref{N=1} is now concluded, thanks to the simple observation $U(t)=t\bM_t(v^E(x, t))$.



\begin{thebibliography}{10}
\bibitem{Alonso2021} Ireneo Peral Alonso and Fernando Soria de Diego, {\em Elliptic and parabolic equations involving the {Hardy}-{Leray} potential}, {De Gruyter Ser. Nonlinear Anal. Appl.}, ISSN: {0941-813X}, Vol. {\bf 38} (2021),  ISBN: {978-3-11-060346-0; 978-3-11-060627-0}.


\bibitem{NLDP}{F. Andreu-Vaillo, J. M. Maz\'{o}n, J.D. Rossi and J. J. Toledo-Melero}, {\em Nonlocal diffusion problems}, {Mathematical Surveys and Monographs},  Vol. {\bf 165}, {American Mathematical Society, Providence, RI; Real Sociedad
              Matem\'{a}tica Espa\~{n}ola, Madrid}, {2010}, {xvi+256}.       
              
\bibitem{Baz} E. G. Bazhlekova, {\em Subordination principle for fractional evolution equations}, Fract.
Calc. Appl. Anal., {\bf 3} (2000), 213--230.

\bibitem{Bern} Ber, J. \emph{Subordinators: Theory and Applications}; Cambridge University Press: Cambridge, UK, 1996.

\bibitem{Ber} Bertoin, J. \emph{L\'{e}vy Processes of Cambridge Tracts in Mathematics}; {Cambridge University Press: Cambridge, UK}, 1996; {Volume 121}.

\bibitem{Bingham} N.H. Bingham, C.M. Goldie and  J.L. Teugels, {\em Regular Variation, Encyclopedia of Mathematics and Its Applications}, Vol. {\bf 27}, Cambridge University Press, Cambridge (1987).

\bibitem{BG1960}{R. M. Blumenthal and R. K. Getoor},
     {\em Some theorems on stable processes}, {Trans. Amer. Math. Soc.}, {\bf 95} (1960), {263--273}.









\bibitem{Feller} W. Feller, {\em An Introduction to Probability Theory and Its Applications}, Vol. II, 2nd edn., Wiley, New York, 1971.


\bibitem{KKS20}{Anatoly N. Kochubei, Yuri G. Kondratiev and Jos{\'e} L. da Silva}, {\em Random time change and related evolution equations. Time asymptotic behavior}, {Stoch. Dyn.}, {\bf 20} (2020) 2050034, 24 Pages.


\bibitem{KKS21}{Anatoly N. Kochubei, Yuri G. Kondratiev and Jos{\'e} L. da Silva}, {\em On fractional heat equation}, {Fract. Calc. Appl. Anal.}, {\bf 24} (2021), {73--87}.

\bibitem{Koch} {Anatoly N. Kochubei}, {\em Distributed order calculus and equations of ultraslow diffusion}, J. Math. Anal. Appl., {\bf 340} (2008) 252--281.



\bibitem{KGW} A. Kumar, J. Gajda and  A. Wyloma\'{n}ska, {\em Fractional Brownian Motion Delayed by Tempered and Inverse Tempered Stable Subordinators}, {Methodol Comput Appl Probab,} \textbf{21}, (2019), 185--202.


\bibitem{Laister}{R. Laister and M. Sier\.{z}ega}, {\em Well-posedness of semilinear heat equations in {$L^1$}}, {Ann. Inst. H. Poincar\'{e} C Anal. Non Lin\'{e}aire}, {\bf 37} (2020), {709--725}.





 \bibitem{MS2012} Mark M. Meerschaert and Alla Sikorskii,  {\em Stochastic Models for Fractional Calculus}, De Gruyter Studies in Mathematics,
    Vol. {\bf 43}, {Second Edition}, {De Gruyter, Berlin},  {2019}.

\bibitem{Mee2013} Mark M. Meerschaert, E. Nane and  P. Vellaisamy, {\em Transient anomalous sub-diffusion on bounded domains}, {Proc. Amer. Math. Soc.}, {\bf 141} ({2013}), {699--710}.

\bibitem{MK2000} R. Metzler and J. Klafter {\em The random walk's guide to anomalous diffusion: a fractional dynamics approach}, Physics Reports
{\bf 339} (2000), 1--77.



\bibitem{Bernstein}{Ren\'{e} L. Schilling, R. Song and Z. Vondra\v{c}ek}, {\em Bernstein functions}, {De Gruyter Studies in Mathematics}, Vol. {\bf 37},{Second Edition}, {Theory and applications}, {Walter de Gruyter \& Co., Berlin}, {2012}.






 


\bibitem{Vaz}{J. L.- V{\'a}zquez}, {\em Asymptotic behaviour for the fractional heat equation in the {Euclidean} space}, {Complex Var. Elliptic Equ.},{\bf 63} (2018), {1216--1231}.



\end{thebibliography}
\end{document}